\font\textmsbm= msbm10 scaled 1200
\font\scriptmsbm= msbm7 scaled 1200
\font\scriptscriptmsbm= msbm5 scaled 1200
\newfam\tmwfama
\textfont\tmwfama\textmsbm
\scriptfont\tmwfama\scriptmsbm
\scriptscriptfont\tmwfama\scriptscriptmsbm
\def\bbb{\fam\tmwfama\textmsbm}
\def\d{\delta}
\def\e{\epsilon}

\def\th{\theta}

\def\l{\lambda}

\def\s{\sigma}

\def\o{\omega}

\def\O{\Omega}
\def\O{\Omega}

\def\pa{\partial}

\def\sm{\setminus}

\def\sb{\subset}      \def\sbe{\subseteq}

\def\tR{\mathchoice{\mbox{\bf R}}{\mbox{\bf R}}{\mbox{\scriptsize \bf
R}}{\mbox{\tiny\bf R}}}
\def\tE{\mathchoice{\mbox{\bf E}}{\mbox{\bf E}}{\mbox{\scriptsize \bf
E}}{\mbox{\tiny\bf E}}}

\def\tC{\mathchoice{\mbox{\bf C}}{\mbox{\bf C}}{\mbox{\scriptsize \bf
C}}{\mbox{\tiny\bf C}}}

\def\tP{\mathchoice{\mbox{\bf P}}{\mbox{\bf P}}{\mbox{\scriptsize \bf
P}}{\mbox{\tiny\bf P}}}

\def\bydef{\,\lower-.1ex\hbox{\rm :}\!=}

\def\cF{{\cal F}}

\def\cJ{{\cal J}}

\def\dist{\mbox{ \rm dist} }
\def\Id{\mbox{ \rm Id} }

\def\proof{\goodbreak\noindent{\sc Proof. }\nobreak}
\def\endproof{\par\nobreak\hbox to \hsize{\hfil\vrule width 5pt height
5pt}\goodbreak\vskip 3pt}
\newtheorem{theor}{Theorem}

\newtheorem{lemma}[theor]{Lemma}

\def\bs{\vskip 0.4cm}

\def\Re{\mbox{ \rm Re }}
\def\Id{\mbox{ \rm Id }}
\def\BMO{\mathchoice{\mbox{ \rm BMO }}{\mbox {\rm BMO}}{\mbox{\scriptsize
\rm BMO}}{\mbox{\tiny \rm BMO}}}
\def\bT{{\bbb T}}
\def\l{\lambda}
\documentstyle[12pt,titlepage,oldlfont]{article}
\makeatletter

\title {Two Remarks on Marcinkiewicz decompositions by Holomorphic Martingales}
\author{Paul F.X. M\"uller\\
Department of Mathematics\\
 J. Kepler Universit\"at\\
A-4040 Linz, AUSTRIA}
\date{October 11, 1993}

\titlepage
\begin{document}
\maketitle
\newpage
\pagestyle{plain}
\pagenumbering{arabic}

\section {Introduction}

The real part of $H^\infty(\bT)$ is not dense in $L^\infty_{\tR}(\bT)$. The
John-Nirenberg theorem in combination with the Helson-Szeg\"o theorem and
the Hunt Muckenhaupt Wheeden theorem has been used to determine whether
$f\in L^\infty_{\tR}(\bT)$ can be approximated by $\Re H^\infty(\bT)$ or
not:
$\dist(f,\Re H^\infty)=0$ if and only if for every $\e>0$ there exists
$\l_0>0$ so that for $\l>\l_0$ and any interval $I\sbe \bT$.
$$|\{x\in I:|\tilde f-(\tilde f)_I|>\l\}|\le |I|e^{-\l/ \e},$$
where $\tilde f$ denotes the Hilbert transform of $f$. See [G] p. 259.
This result is contrasted by the following

\begin{theor} Let $f\in L^\infty_{\tR}$ and $\e>0$. Then
there is a function $g\in H^\infty(\bT)$ and a set $E\sb \bT$ so that
$|\bT\sm E|<\e$ and $$f=\Re g\quad\mbox{ on } E.$$  \end{theor}

This theorem is best regarded as a corollary to Men'shov's correction theorem.
For the classical proof of Men'shov's theorem see [Ba, Ch VI \S 1-\S4].

Simple proofs of Men'shov's theorem -- together with significant extensions --
have been obtained by S.V. Khruschev in [Kh]
and S.V. Kislyakov in [K1], [K2] and [K3].

In [S] C. Sundberg used $\bar\pa$-techniques (in particular [G, Theorem  VIII.1.
gave a proof of Theorem 1 that does not mention Men'shov's theorem.

The purpose of this paper is to use a Marcinkiewicz decomposition on Holomorphic
Martingales to give another proof of Theorem 1. In this way we avoid uniformly
convergent Fourier series  as well as $\bar\pa$-techniques.

Holomorphic Martingales enter in the proof of the following lemma.

\begin{lemma} There exist $c_1,c_2>0$ so that for every $f\in \BMO A$, where
$||f||\le 1$, $\e>0$ and $\l\in \tR^+$ there exists $g\in H^\infty(\bT)$
and $E\sb \bT$
\begin{eqnarray*}||g||_\infty&\le& \l \\
|f(\th)-g(\th)|&\le&\e\quad\mbox{ on } E\\
|\bT\sm E|&\le&\frac 1\e e^{-\l c_1} c_2.\end{eqnarray*}
\end{lemma}

Consider complex Brownian Motion $(z_t)_{t\ge 0}$ on the Wiener space
$(\O,(\cF_t),\cF,\tP)$. A complex valued random variable $X$ on $\O$ is
called holomorphic if the conditional expectation $$X_t=\tE(X|\cF_t)$$ admit
a stochastic integral representation of the form $X_t=X_0+\int^t_0F_s dz_s$,
where $F_s$ is  adapted to $\cF_s$.

$H^p(\O)$ denotes the closure in $L^p(\O)$ of holomorphic random variables.
$\BMO(\O)$ denotes the closure of holomorphic random variables under the norm
$$\sup_t||\tE(|X-X_t| \, |\cF_t)||_\infty.$$

The connection to analytic functions is provided by operators $M,N$ so that
$$\begin{array}{ccccc}
H^p(\bT)&&\stackrel {\Id}{\longrightarrow}&&H^p(\bT)\\
M&\searrow&&\nearrow&N\\
&&H^p(\O)&&\end{array}$$
where $||M||_p=||N||_p=1$ and
$$\begin{array}{ccccc}
\BMO A(\bT)&&\stackrel {\Id}{\longrightarrow}&&\BMO A(\bT)\\
M&\searrow&& \nearrow &N\\
&&\BMO (\O)&&\end{array}$$
where $||M||_{\BMO}\le C_0, ||N||_{\BMO}\le C_0$.

These probabilistic ideas have a quite long history and were useful in
several problems of Analysis. See [F], [G-S], [Ma]
and [V].

\section{Proofs of the results}

{\bf Proof of Lemma 2.}
Fix $\l>0$ and let
\begin{eqnarray*}   \s&=&\inf \{t:|z_t|>1\}\\
\tau&=&\inf\{t\le\s:|f(z_t)|>\l\}\\
G_t&=& f(z_{t\wedge \tau})\\
F_t&=& f(z_t)\\
g(\th)&=& N(G)(\th).\end{eqnarray*}
Then \begin{eqnarray*} ||g||_\infty&\le&\l\\
g&\in& H^\infty\end{eqnarray*}
and \begin{eqnarray*} |\{\th:|f-g|>\e\}|&\le&\frac 1\e||f-g||_1\\
&\le&\frac 1\e||N(F-G)||_1.\end{eqnarray*}
By [M, Lemma 1] we get $$||N(F-G)||_1\le||F-G||_1\le 2\int_{\{F^*>\l\}}|F|dP$$
where $F^*=\sup|F_t|$. By Cauchy-Schwartz we obtain

$$\int_{\{F^*>\l\}}|F|dP\le\{F^*\ge\l\}^{1/2}c_2\le e^{-\l c_1}c_2,$$
because $F\in \BMO(\O)$ implies $\tE e^{F^*c}<\infty$. This implies the
estimate
$$|\{\th:|g-f|>\e\}|\le\frac 1\e e^{-\l c_1}c_2.$$
\bs

{\bf Proof of Theorem 1.}
Given $\e>0$ we select $\l_n\in\tR^+$ so that $$\sum^\infty
_{n=0}e^{-\l_nc_1}c_22^n<\e$$
and $$\sum^\infty_{n=0}\l_n 2^{-n}<\infty.$$
Then given a function $h:\bT\to \tC$ and $\d>0$ we define
$$T_\d(h)(\th)=\left\{\begin{array}{ccc} h(\th)&\mbox{if}&|h(\th)|\le\d\\
\d&\mbox{if}&|h(\th)|\ge\d\end{array}\right.$$
Now consider $u_0\in L^\infty_{\tR}(\bT)$ with $||u_0||_\infty=1$ and let
$\tilde u_0$ be the Hilbert transform of $u_0$ then $u_0+i\tilde u_0\in \BMO A$
and
$$||u_0+i\tilde u_0||_{\BMO A}\le C||u_0||_\infty.$$
We next apply an interation procedure from [S].
\bs
{\bf Step 1.} Use Lemma 2  to obtain $E_1\sb \bT$, $g_1\in H^\infty$
with $||g_1||_\infty\le\l_1$ so that $$|u_0+i\tilde
u_0-g_1|<1/2\quad\mbox{on } E_1$$ and $$|\bT\sm E_1|\le 2e^{-\l_1c_1}c_2.$$

{\bf Induction Step.}\,  We have already constructed $u_0,\dots,u_{n-1}\in
L^\infty_{\tR}, g_1,\dots,g_n\in H^\infty(\bT)$ and $E_1,\dots, E_n\le \bT$
so that for $j\le n$
$$||g_j||_\infty\le\l_j2^{-j}$$
$$|u_{j-1}+i\tilde u_{j-1}-g_j|\le 2^{-j}\quad\mbox{on } E_j$$
$$|\bT\sm E_j|\le e^{-c_1\l_1}c_22^j.$$
Now we let $$u_{n}:=T_{2^{-n}}(u_{n-1}-\Re g_n)$$
and we have $$u_n= u_{n-1}-\Re g_n\quad\mbox{on }E_n$$
$$||u_n||_\infty\le 2^{-n}.$$
By Lemma 2 we find $g_{n+1}\in H^\infty(D)$, $E_{n+1}\sbe\bT$ so that
$$||g_{n+1}||_\infty\le\l_{n+1}2^{-n-1},$$
$$|u_n+i\tilde u_n-g_{n-1}|\le 2^{-n-1}\quad\mbox{on } E_{n+1},$$
$$|\bT\sm E_{n+1}|\le e^{-c_1\l_{n+1}}c_22^{n+1}.$$

Having completed the construction we set $$g:=\sum^\infty_{j=1}g_j$$ which
defines an element in $H^\infty(\bT)$.
Tracing back we see that
$$\sum^\infty_{n=1}u_n=\sum^\infty_{n=0}u_n-\sum^\infty_{n=1}\Re
g_n\quad\mbox{on } \bigcap^\infty_{n=1}E_n$$
or $$u_0=\Re g\quad\mbox{on } \bigcap^\infty_{n=1}E_n.$$
It remains to estimate $|\bigcap^\infty_{n=1}E_n|$ from below:

\begin{eqnarray*}
\biggl\vert\bigcap^\infty_{n=1} E_n\biggr\vert&\ge&|\bT|-\sum_{n=1}|\bT\sm
E_n|\\
&\ge&|\bT|-\sum^\infty_{n=1}e^{-\l_nc_1}2^nc_2\\
&\ge&|\bT|-\e.\end{eqnarray*}

\bs

\section{A Refinement of Lemma 2}

In the above argument we gave just an estimate for the {\bf size} of the
set $$\{\th:|f(\th)-g(\th)|<\e\}$$
but did not give any indication where to find this set.
A more detailed analysis of the \lq\lq conditional expectation" operator
$N$ gives estimates which relate the probabilistic Marcinkiewicz
decomposition to classical maximal functions.
\bs
Let $h:\bT\to\tC$ be a function, then let $h^\#$ be the non tangential
maximal function and define
$$M_{HL}(h)(\th):=\sup_I\int_I|h|\frac{dt}{|I|}$$ where the $\sup$ is taken over
intevals in $\bT$ which contain $\th$.
 Let $g$ be defined as in the proof of Lemma 2 then we have the pointwise estima

\begin{theor} 1. $|f(\th)-g(\th)|\le C(|f(\th)|+\l)M_{HL}(\chi_{H_\l})(\th)$ whe
$H_\l=\{f^\#>\l\}$.

2. Let $f\in BMO$, with $||f||\le 1$, then
for every $N>0$ there exists $\l>0$ and $B\sb\{\th\in\bT:|f(\th)|\le N\}$
so that $$|\bT\sm B|\le e^{-Nc_1}c_2,$$
$$M_{HL}(\chi_{H_\l})(\th)\le e^{-\l c_3}\quad\mbox{ for } \th\in B.$$
\end{theor}

\proof ad 1.
For $\th\in\bT$ and $z\in D$ let
$$P_\th(z):=\frac{1-|z|^2}{|e^{i\th}-z|^2}.$$
Fix $0<r< 1$ consider the stopping times $$\s_r:=\inf \{t:|z_t|>r\}$$ and
let $$F_\l:=\{\o\in\O:\tau(\o)<\s(\o)\}$$   $$E_\l:=\{z\in D:|f(z)|>\l\}$$
where the stopping time $\tau$ has been defined during the proof of Lemma
2. Then  for $\th\in\bT$ we have (using formula (1) in [Du, Section 3,2])
\begin{eqnarray*} g(\th)&=&N(F-G)(\th)\\
&=&\lim_{r\to 1}\tE((f(z_{\s_r})-f(z_{\s_r\wedge \tau}))P_\th(z_r))\\
&=&\lim_{r\to
1}\tE((f(z_{\s_r})-f(z_{\s_r\wedge\tau}))\chi_{F_\l}P_\th(z_r))\\
&\le&(|f(\th)|+\l)\lim_{r\to 1}\tE(\chi_{F_\l}P_\th(z_{\s_r})).\end{eqnarray*}
For $A\sb D$ let $$\o(A):=\tP\{z_t\in A,\mbox{ for some } t<\s\}.$$
Then  $\o(E_\l)=\tP(F_\l)$ and by the strong Markov Property we have: (see
[D],p. 923 or [V], p. 112) $$\lim_{r\to
1}\tE(\chi_{F_\l}P_\th(z_{\s_r}))=\int_{\pa E_\l}P_\th(z)\o(dz).$$
The integral on the RHS is called balyage or sweep of $\o_{|\pa E_\l}$ and
has been much studied because of its relation to Carleson-measures and
$\BMO$.
See [G], pp. 229, 239 and 240.
The argument in [G], p. 239 gives the estimate $$\int_{\pa
E_\l}P_\th(z)\o(dz)\le C_3\sup_{0\le h<1}\frac{\o(\pa E_\l\cap S_h)}{h}$$
where $$S_h:=\{re^{i\psi}:1-h\le r< 1,|\psi-\th|\le h\}.$$
The result of Burkholder Gundy Siverstein gives for every harmonic
function $u:D\to \tR$ $$\o\{u>\l\}\le C|\{u^\#>\l\}|.$$
See [P], p. 36.
Therefore by [G, Lemma I.5.5] $\o$ is a Carleson Measure.

Hence a simple stopping time argument gives for every $0\le h\le 1$
$$\o(\pa E_\l\cap S_h)\le C|\{f^\#>\l\}\cap 3I_h|$$ where $I_h=S_h\cap\bT$.
We therefore have the estimate.

$$\int_{\pa E_\l}P_\th(z)\o(dz) \le C\sup_{0\le h\le
1}\frac{|\{f^\#>\l\}\cap 2I_h|}{h}$$
And by choice of $I_h$ the LHS is dominated by $$CM_{HL}(\chi_{H_\l})(\th)$$ whe
$H_\l=\{\psi\in\bT:f^\#(\psi)>\l\}$.\bs
ad 2. As $f\in \BMO$ there exists $\d_0>0$ and $C_0>0$ so that for each
$\l>0$
 $$|\{f^\#>\l\}|\le e^{-\l\d_0}C_0.$$
Now choose $\d_1=\d_0/2$ and let
\begin{eqnarray*} H&=&\{f^\#>\l\}\\
G&=&\{|f|<N\}\\
\cJ&=&\{I\sb\bT:|H\cap I|>e^{-\l \d_1}|I|, \, I\mbox{ Intervall}\}\\
J&=& \bigcup_{I\in\cJ}I\\
B&=&G\sm J.\end{eqnarray*}
The weak type $1:1$ bound for the Hardy Littlewood maximal function gives
$$|J|\le |H|e^{\l\d_1}C\le e^{-\l\d_1}C.$$

Hence \begin{eqnarray*} |\bT\sm B|&\le&|\bT\sm G|+|J|\\
&\le& e^{-N\d_0}C_0+e^{-\l\d_1}C.\end{eqnarray*}

Moreover, by definition, we have $$M_{HL}(\chi_H)(\th)\le e^{-\l\d_1}\mbox{
for } \th\in B$$ and this completes the proof.
\endproof


\newpage
\begin{thebibliography}{B-G-S}
\bibitem[Ba]{} Bary N.K., A Treatise on Trigonometric Series, Pergamon Press
(1964).

\bibitem[B-G-S]{} Burkholder D.L., Gundy R.F., Silverstein M.L., A maximal
function Charactersation of the class $H^p$, Trans. Amer. Math. Soc. 157
(1971), 137--153.

\bibitem[D]{} Davis B., Brownian Motion and analytic Functions, Annals. of
Prob. 7 (1979), 913--932.

\bibitem[Du]{} Durrett R., Brownian Motion and Martingal in Analysis, The
Wadsworth Mathematical Series (1984).

\bibitem[F]{} F\"ollmer H., Stochastic Holomophy, Math. Ann. 207 (1974),
245--255.

\bibitem[G]{} Garnett J.B., Bounded analytic Functions, Academic Press,
(1981).

\bibitem [G-S]{} Getoor R.K., Sharpe M.J., Conformal Martingales, Inv. Math.
16 (1972), 271--308.

\bibitem[Kh]{} Khrushchev S.V., Men'shov's correction theorem and Gaussian
processes, Proc. of Steklov Institute of Math. 155 (1983), 147--175, (Engl.
Transl.).

\bibitem[K1]{} Kislyakov S.V., The Fourier Coefficients of the boundary Values
of Functions analytic in the Disk and in the Bidisk, Proc. of  Steklov
Institute of Math. 155 (1983), 75--91, (Engl. Transl.).

\bibitem[K2]{} Kislyakov S.V., Quantitative aspects of correction theorems, Zap.
Naucn. Sem. Leningrad, Otdel. Mat. Inst. Steklov (LOMI) 92 (1979), 182--191.

\bibitem[K3]{}  Kislyakov S.V., A new correction theorem, Math. USSR Izvestiya
24 (1985) 283--305 (Engl. Transl.).

\bibitem[M]{} M\"uller P.F.X., Holomorphic Martingales and Interpolation
between Hardy Spaces, appears in Jour. d'Analyse Math. (Jerusalem).

\bibitem[Ma]{} Maurey B., Isomorphismes entre espaces $H^1$, Acta Math. 145
(1980), 79--120.

\bibitem[P]{} Petersen K.E., Brownian Motion, Hardy Spaces and Bounded mean
Oscillation, L.M.S. Lecture Note Series, 28, Cambridge University Press
(1977).

\bibitem[S]{} Sundberg C., Truncation of $\BMO$ Functions, Ind. Univ. Math.
J. 33 (1984), 749--771.

\bibitem[V]{} Varopoulos N.Th., The Helson Szeg\"o Theorem and
$A_p$-Functions for Brownian Motion and Several Variables, Jour. of Funct.
Anal. 39 (1980), 85--121.\end{thebibliography}
\end{document}